\documentclass[12pt]{article}
\usepackage{e-jc}

\usepackage{amsmath}

\dateline{May 12, 2018}{TBD}{TBD}
\Copyright{The authors. Released under the CC BY-ND license (International 4.0).}

\title{Chromatic Vertex Folkman Numbers\thanks{
Supported by the National Natural Science Foundation (11361008)
and the Guangxi Natural Science Foundation (2011GXNSFA018142).}
}

\author{Xiaodong Xu\\
\small Guangxi Academy of Sciences\\[-0.8ex]
\small Nanning, 530007, P.R.China\\[-0.8ex]
\small\tt xxdmaths@sina.com\\
\and
Meilian Liang\\
\small School of Mathematics and Information Science\\[-0.8ex]
\small Guangxi University, Nanning, 530004, P.R.China\\[-0.8ex]
\small\tt gxulml@163.com\\
\and
Stanis{\l}aw P. Radziszowski\\
\small Department of Computer Science\\[-0.8ex]
\small Rochester Institute of Technology\\[-0.8ex]
\small Rochester, NY 14623, USA\\[-0.8ex]
\small\tt spr@cs.rit.edu}

\begin{document}
\maketitle
\begin{abstract}
For graph $G$ and integers
$a_1 \ge \cdots \ge a_r \ge 2$, we write
$G \rightarrow (a_1 ,\cdots ,a_r)^v$  if and only if for every
$r$-coloring of the vertex set $V(G)$ there exists a monochromatic
$K_{a_i}$ in $G$ for some color $i \in \{1, \cdots, r\}$.
The vertex Folkman number $F_v(a_1 ,\cdots ,a_r; s)$
is defined as the smallest integer $n$ for which there
exists a $K_s$-free graph $G$ of order $n$ such that
$G \rightarrow (a_1 ,\cdots ,a_r)^v$. It is well known that
if $G \rightarrow (a_1 ,\cdots ,a_r)^v$ then $\chi(G) \geq m$,
where $m = 1+ \sum_{i=1}^r (a_i - 1)$.
In this paper we study such Folkman graphs $G$ with
chromatic number $\chi(G)=m$, which leads to
a new concept of chromatic Folkman numbers.
We prove constructively some existential results,
among others that
for all $r,s \ge 2$ there exist $K_{s+1}$-free graphs $G$
such that $G \rightarrow (s,\cdots_r,s)^v$ and $G$ has
the smallest possible chromatic number $r(s-1)+1$
for this $r$-color arrowing to hold.
We also conjecture that, in some cases, our construction is
the best possible, in particular that for every $s \ge 2$ there exists
a $K_{s+1}$-free graph $G$ on $F_v(s,s; s+1)$ vertices
with $\chi(G)=2s-1$ such that $G \rightarrow (s,s)^v$.
\end{abstract}


\section{Preliminaries, notation and definitions}
\label{sintro}
Throughout this paper, we consider only
finite undirected loopless simple graphs.
For graph $G=(V,E)$, denote by $V(G)$ the set of its
vertices, and by $E(G)$ the set of its edges.
A complete graph of order $n$ is
denoted by $K_n$, and a cycle of length $n$ by $C_n$.
The clique number of $G$ is denoted by $cl(G)$,
and the chromatic number by $\chi(G)$. An $(s,t)$-graph
is a graph that contains neither an $s$-clique
nor a $t$-independent set.

The set $\{1, \cdots, n\}$ is denoted by $[n]$.
Let $r, s, a_1, \cdots, a_r$ be integers such
that $r \ge 2$,
$s >  \max \{a_1, \cdots, a_r \}$ and
$\min \{a_1, \cdots, a_r \} \ge 2$.
We write $G \rightarrow (a_1 ,\cdots ,a_r)^v$
if and only if for every
$r$-coloring of $V(G)$ there exists
a monochromatic $K_{a_i}$ in $G$ for some color
$i \in [r]$.
The sets of vertex Folkman graphs are defined as
$$\mathcal{F}_v(a_1, \cdots, a_r; s)=\{G \;|\; G \rightarrow
(a_1, \cdots, a_r)^v \textrm{ and } cl(G)<s\},$$
and the vertex Folkman numbers
are defined as the smallest orders of graphs in these sets,
namely
$$F_v(a_1, \cdots, a_r;s)=\min\{|V(G)|\;|\;G\in
\mathcal{F}_v(a_1, \cdots, a_r;s) \}.$$

\medskip
In 1970,
Folkman \cite{Folkman} proved that for every
$s > \max \{a_1, \cdots, a_r\}$ the sets\break
$\mathcal{F}_v(a_1, \cdots, a_r; s)$ are nonempty.
If $a_1 = \cdots = a_r =a$, then we will use a simpler
notation for the corresponding set of vertex Folkman graphs
$\mathcal{F}(r,a,s) = \mathcal{F}_v(a_1 ,\cdots ,a_r; s)$,
and for numbers $F(r,a,s) = F_v(a_1 ,\cdots ,a_r; s)$.
The case of $F(r,s,s+1)$ was studied in particular
by Dudek and R\"{o}dl \cite{DudekRodl2010jctb} and
H\`{a}n, R\"{o}dl and Szab\'{o} \cite{HRS}.
The latter work contains the result stated in
the following theorem.

\medskip
\begin{theorem} \label{pupperbound}
{\rm (H\`{a}n-R\"{o}dl-Szab\'{o}, 2018)}
For any positive integer $r$ there exists a constant $C = C(r)$
such that for every $s \geq 2$ it holds that $F(r,s,s+1) \leq Cs^2 \log^2s$.
\end{theorem}

Set $m = 1+ \sum_{i=1}^r (a_i - 1)$. It is well known that
if $G \rightarrow (a_1, \cdots, a_r)^v$ then
$\chi(G) \geq m$ \cite{Nenov04}.
In this paper we study vertex Folkman graphs $G$ with
the corresponding Folkman numbers, when the graphs
$G$ satisfy an additional constraint on their chromatic number,
namely $\chi(G)=m$.
This motivates the concept of the {\em minimum chromatic Folkman
graphs}, $\mathcal{F}_v^{\chi}(a_1, \cdots, a_r; s)$,
and corresponding {\em minimum chromatic Folkman numbers},
$F_v^{\chi}(a_1 ,\cdots ,a_r;s)$. Formally:
$$\mathcal{F}_v^{\chi}(a_1, \cdots, a_r; s)=
\{G \;| \; G \in \mathcal{F}_v(a_1, \cdots, a_r; s) \textrm{ and } \chi(G) = m\},$$
and
$$F_v^{\chi}(a_1 ,\cdots ,a_r;s) =
\min\{|V(G)|\;| \; G\in \mathcal{F}_v^{\chi}(a_1, \cdots, a_r;s) \}.$$

\medskip
The remainder of this paper focuses on minimum chromatic
vertex Folkman graphs and numbers, and we will call them simply
chromatic Folkman graphs and numbers, respectively.
Similarly as in the classical case, we will
use further notational abbreviation for
the diagonal chromatic cases,
namely, if $a_1 = \cdots = a_r =a$, then we set
$\mathcal{F}^{\chi}(r,a,s) = \mathcal{F}_v^{\chi}(a_1 ,\cdots ,a_r;s)$
and $F^{\chi}(r,a,s) = F_v^{\chi}(a_1 ,\cdots ,a_r;s)$.

\bigskip
The main result of this paper, Theorem 3 in Section 2.1,
shows the existence of $F^{\chi}(r,s,s+1)$
for $r,s \ge 2$ by an explicit construction.
We must admit that the implied upper bound is much larger than
the bound in Theorem 1, but our construction gives an additional
constraint on the chromatic number.
We also conjecture that, in some cases, our
construction is the best possible one, in particular that for
every $s \ge 2$ there exists a $K_{s+1}$-free graph $G$ on
$F_v(s,s; s+1)$ vertices with chromatic number
$\chi(G)$ equal to $m=2s-1$
such that $G \rightarrow (s, s)^v$, or equivalently,
that $F_v^{\chi}(s,s;s+1)=F_v(s,s;s+1)$.
We wish to remark that a technique as in
this paper should lead to more general existence results
for $F_v^{\chi}(a_1 ,\cdots ,a_r;s)$, and with
some further enhancements also for generalized Folkman numbers,
where one avoids monochromatic graphs $G_i$ instead of $K_{a_i}$.
However, we do not study such extensions in this work.

\medskip
In another direction, one can consider coloring the edges
instead of vertices, which leads to concepts analogous
to those above for vertices. Let
$M = R(a_1, \cdots, a_r) =
\min\{n~|~K_n \rightarrow (a_1, \cdots, a_r)^e\}$,
i.e. set $M$ to be the classical multicolor
Ramsey number for fixed $a_1, \dots, a_r$
(In this paragraph the subscript/superscript $e$ means that we
color edges instead of vertices).
A classical argument easily proves that
if $G \rightarrow (a_1, \cdots, a_r)^e$
then $\chi(G) \geq M$ \cite{Lin1972a}.
Further, we can define
$\mathcal{F}_e^{\chi}(a_1, \cdots, a_r; s)=
\{G\;| \; G \in \mathcal{F}_e(a_1, \cdots, a_r; s) \textrm{ and } \chi(G) = M\},$
and
$F_e^{\chi}(a_1 ,\cdots ,a_r;s) =
\min\{|V(G)|\;| \; G\in \mathcal{F}_e^{\chi}(a_1, \cdots, a_r;s) \}.$
From the construction by Ne\v{s}et\v{r}il and R\"{o}dl in
\cite{NesetrilRodl} we can see that there exists
a $K_{s+1}$-free graph $G$ with chromatic number
equal to $R(s,s)$ such that $G \rightarrow (s,s)^e$.
For $3 \leq k \leq l$, essentially the same
reasoning as in \cite{NesetrilRodl}
implies the existence of $K_{s+1}$-free
graphs $G$ with $\chi(G) = R(k,l)$
and $G \rightarrow (k,l)^e$.
In a similar way, one could also show that
$F_e^{\chi}(a_1 ,\cdots ,a_r;s)$ exists.

\medskip
In the sequel we will be coloring only vertices.
The rest of this paper is organized as follows.
In Section \ref{Newtype}, the existence of the
chromatic Folkman numbers $F^{\chi}(r,s,s+1)$
is proved, the corresponding upper bound is discussed,
and a conjecture is posed.
In Section \ref{Minimal}, we study some
minimal Folkman graphs with the smallest minimum
degree. Though not yet directly connected to chromatic
Folkman numbers, we point to a possible connection
between them.

\medskip
\section{Chromatic Folkman numbers} \label{Newtype}

Our main motivation to study chromatic Folkman graphs
and numbers is to understand how they may differ from
the regular Folkman cases. If for a special family
of cases we find that the chromatic cases yield
the same numbers as the regular ones, then we could
see it as a stronger version of the Folkman theorem.
In any case, we may be able to see better the structure
of extremal Folkman graphs.

Before we prove the existence of $F^{\chi}(r,s,s+1)$,
we want to observe that $F_v^{\chi}(a_1 ,\cdots ,a_r;s)$ and
$F_v(a_1 ,\cdots ,a_r;s)$ can be different.
It is known that there exists exactly one $K_5$-free
graph $Q$ of order 13 such that $Q \rightarrow (3,4)^v$
\cite{Nenov01}. This graph $Q$ is the complement
of the unique $(3,5)$-Ramsey critical graph,
which is cyclic on the set $\mathbb{Z}_{13}$ with arcs
of length 1 and 5. One can easily check that $\chi(Q)=7$.
Thus, since in this case $m=6$, we have
$F_v^{\chi}(3,4;5)> 13 = F_v(3,4;5)$.
Bikov and Nenov \cite{BNc17} provided
us with other examples
of this type, in particular they observed that
$F_v^{\chi}(4,4;6) > F_v(4,4;6)$, since the
only extremal graph for $F_v(4,4;6)$ is $Q+K_1$.
More such examples follow from their recent work \cite{BN15,BNa17}.
On the other hand we feel that these examples are special
in that they exploit larger difference between arrowed
and avoided graphs. This is captured in the following
problem for the borderline diagonal cases,
for two and more colors.

\bigskip
\noindent
{\bf Problem 1.} {\it
Is it true that $F^{\chi}(r,s,s+1) = F(r,s,s+1)$
for all $r,s \ge 2$}?

\bigskip
Lemma 2 in the next section gives a positive answer to
Problem 1 for all $r \ge 2$ with $s=2$.
A graph presented in Figure 1 in \cite{PRU99} and the
computations described therein implicitly give a positive
answer in one more special 2-color case,
namely it holds that $F^{\chi}_v(3,3;4) = F_v(3,3;4)=14$.
One can ask similar questions related to minimum
chromatic edge Folkman numbers
(when coloring edges instead of vertices),
but these seem much more difficult to answer.

\medskip
While several general cases of vertex Folkman
numbers have been studied, the problem of finding
their exact values for small cases remains elusive.
This holds
even just for two colors when $s$ is smaller than $m-1$.
For instance, the case of $F_v(4,4;5)$
seems to be difficult, for which only the bounds
$19 \leq F_v(4, 4; 5) \leq 23$ are known \cite{BNb17,Fv445}.
The computational approach is often too expensive.
Just testing a single instance of arrowing for
an upper bound witness graph is not easy in most cases,
and improving lower bounds is much harder since it
may involve a very large number of arrowing
instances. Both upper and lower bounds for
edge Folkman numbers
tend to be computationally still harder.

\medskip
\subsection{Main theorem}

\medskip
The upper bound on $F^{\chi}(r,s,s+1)$, which can be obtained
by the construction of Theorem~\ref{largercase} below, is
rather large. This and other bounds we could derive for
$F^{\chi}(r,a,s)$ also seem less tight than those
known for $F(r,a,s)$, except the special case
for $a=2$ and $s=3$
captured by Lemma \ref{Fr23}. We follow this lemma
by a theorem describing the quite special but more
general and more difficult case
of $F^{\chi}(r, s, s+1)$.

\medskip
\begin{lemma}\label{Fr23}
For all $r \ge 2$, $F^{\chi}(r,2,3)$ exists and
it is equal to $F(r,2,3)$.
\end{lemma}

\begin{proof}
The sets $\mathcal{F}^{\chi}(r,2,3)$ consist of triangle-free
graphs with chromatic number equal to $r+1$,
which in the case of arrowing $K_2$ clearly coincides
with the smallest order graphs in $\mathcal{F}_v(r,2,3)$.
The latter sets are known to be nonempty for all $r \ge 2$.
Thus we also have $F^{\chi}(r,2,3) = F(r,2,3)$.
\end{proof}

\bigskip
In the basic case of two colors, one
can easily see that $F^{\chi}(2,2,3)=5$,
because $C_5 \in \mathcal{F}_v(2,2;3)$
and $\chi(C_5)=m=3$. Clearly, we also have
$F_v(2,2;3)=5$. We know that $F_v(2,2,2;3)=11$,
or the smallest 4-chromatic triangle-free graph
has 11 vertices, witnessed by the Gr\"{o}tzsch
graph. Also, it is known that
$F^{\chi}(4,2,3)=F_v(2,2,2,2;3)=22$,
or the smallest 5-chromatic
triangle-free graph has 22 vertices \cite{JR}.
The best known bounds in the first open case,
$32 \le F^{\chi}(5,2,3) \le 40$, are due
to Goedgebeur \cite{JG}.

\medskip
\begin{theorem}\label{largercase}
For given integers $r \ge 2$ and $s \geq 3$,
let $b_i = i(s-1)+1$ for $i \in [r-1]$,
and $B=\prod_{i=1}^{r-1} b_i$.
Then $F^{\chi}(r, s, s+1)$ exists and
$$F^{\chi}(r, s, s+1) \leq 1+s+ \sum_{i=2}^{r-1} F^{\chi}(i,s,s+1) +
B \cdot F^{\chi}(r,s-1,s).\eqno{(1)}$$
In particular, for all $s \geq 3$,
the chromatic Folkman number $F^{\chi}(2, s, s+1)$ exists and
we have $F^{\chi}(2, s, s+1) \leq 1 + s + sF^{\chi}(2, s-1 , s).$
\end{theorem}

\begin{proof}
For $r \ge 3$, we will construct a graph
$G(r,s) \in \mathcal{F}^{\chi}(r,s,s+1)$ given any graphs in
each of $\mathcal{F}^{\chi}(i,s,s+1)$ for $2 \le i \le r-1$
and any graph in $\mathcal{F}^{\chi}(r,s-1,s)$. The vertices
of the graph $G(r,s)$ will be formed by vertices of given graphs
corresponding to the terms of the right-hand-side of (1).
The proof is using simultaneous induction on $r$ and $s$,
and it has two main parts: construction of $G(r,s)$, and
the proof that $G(r,s)$ has required properties.

Note that the second part of the theorem is just
an instantiation of the first part for two colors,
$r=2$, in which case the main summation of (1) is empty.
Thus, the basis of our induction is formed by the sets
$\mathcal{F}^{\chi}(i,2,3)$, which are nonempty by
Lemma \ref{Fr23}, and where the corresponding Folkman
numbers satisfy $F^\chi(i,2,3)=F(i,2,3)$ for all $i \ge 2$.

\bigskip
\noindent
{\bf Construction of the graph $G(r,s)$ by induction for $r,s \ge 3$.}

\medskip
Let $G_0$ be the graph of order 1, $G_1 = K_s$,
and set $V_0=V(G_0)$, $V_1=V(G_1)$.
We may assume that the graphs $G_i=(V_i,E_i)$
such that $|V_i| = F^{\chi}(i,s,s+1)$ and
$G_i \in \mathcal{F}^{\chi}(i,s,s+1)$
have been already constructed, for $2 \le i < r$.
Therefore, we know that $\chi(G_i) = b_i = i(s-1)+1$
for $i \in \{0, \cdots, r-1\}$, and hence we can partition
each of the sets of vertices $V_i$ into
$\chi(G_i)$ nonempty independent sets
$V_i(j)$ in $G_i$, so that
$$V_i = \bigcup_{j=1} ^{\chi(G_i)} V_i(j).$$

Let $H$ be any graph in the set
$\mathcal{F}^{\chi}(r,s-1,s)$ with
$|V(H)|=F^{\chi}(r,s-1,s)$ vertices.
For the part of $G(r,s)$ corresponding to the last term of (1),
we take $B=\prod_{i=1}^{r-1} b_i$ isomorphic copies
$H(j_0, \cdots, j_{r-1})$ of $H$
indexed by $r$-tuples $(j_0, \cdots, j_{r-1})$, where
$1 \leq j_k \leq \chi(G_k)$ for $k \in \{0, \cdots, r-1\}$.
Note that the order of $G(r,s)$ is equal to the
right-hand-side of (1) because the sizes of parts described
above match exactly its terms,
$$V=V(G(r,s)) = V_0 \;\cup\; V_1 \;\cup\;
\bigcup_{i=2} ^{r-1} V_i \;\cup
\bigcup_{(j_0, \cdots, j_{r-1})}
V(H(j_0, \cdots, j_{r-1})).$$

Finally, we complete the construction of $G(r,s)$
by adding the edges with one end in any
of the sets of vertices $V(H(j_0, \cdots, j_{r-1}))$
and the other end in the sets $V_i(j)$, as follows:
For each fixed $r$-tuple $(j_0, \cdots, j_{r-1})$, where
$1 \le j_i \le b_i$ for $i \in \{0, \cdots, r-1\}$,
we add all possible edges with one end in
$V(H(j_0, \cdots, j_{r-1}))$ and the other end in $V_i(j_i)$.

\bigskip
\noindent
{\bf Proof that $G(r,s) \in \mathcal{F}^{\chi}(r,s,s+1)$.}

\medskip
We need to show that for all $r,s \ge 3$ we have:
(i) $cl(G(r,s))<s+1$,
(ii) in every $r$-coloring
of the vertices $V$ we have a monochromatic $K_s$,
and (iii) $\chi(G(r,s))=m=r(s-1)+1$.

\medskip
\noindent
(i)
Assume contrary, and suppose that some set $S\subset V$
of order $s+1$ induces $K_{s+1}$. Let $k_i=|S \cap V_i|$
for $0 \le i <r$. From the construction we see that $k_i\le s$
and there exists exactly one $t$ for which $k_t>0$.
Similarly, there exists exactly one $r$-tuple
$(j_0, \cdots, j_{r-1})$ such that
$h = |S \cap V(H(j_0, \cdots, j_{r-1}))| >0$.
Note that since $h \le s-1$ and $k_t+h=s+1$,
then $k_t \ge 2$. However, each vertex of
$H(j_0, \cdots, j_{r-1})$ can be adjacent only
to independent sets $V_i(j) \subset V_i$,
hence we have a contradiction for $i=t$, and
thus $cl(G(r,s)) \le s$. We observe that actually
$cl(G(r,s))=s$ since the graph $G(r,s)$
contains $G_1=K_s$.

\medskip
\noindent
(ii)
Assume contrary, and suppose that some $r$-coloring
$C$ of $V$ does not contain any monochromatic $K_s$.
Let $c_0$ be the color of $V_0$, and $c_1$ be a
different color of one of the vertices in $V_1$.
Using the assumptions that
$G_i \in \mathcal{F}^{\chi}(i,s,s+1)$
and that there is
no monochromatic $K_s$ in $C$ restricted to $V_i$,
we can see that the $(i+1)$-st color $c_i$
must be used for some vertex in $V_i$.
Thus, we can find a vertex $v_i \in V_i$ in a new color
$c_i=C(v_i)$, for each $2 \le i <r$.
Let $j_i'$ be such that $v_i \in V_i(j_i')$.
The graph $H(j_0', \cdots, j_{r-1}')$
is isomorphic to $H \in \mathcal{F}^{\chi}(r,s-1,s)$,
and therefore it contains a monochromatic $(s-1)$-clique
$S$ in color $c_i$ for some $0 \le i <r$.
Now in $G(r,s)$, this $S$ can be extended to
a monochromatic $K_s$ in color $c_i$ by adding
vertex $v_i$, which is a contradiction.

\medskip
\noindent
(iii)
Part (ii) implies that $\chi(G(r,s)) \geq r(s-1)+1$,
hence we only need to prove that
$\chi(G(r,s)) \leq r(s-1)+1$.
We will show how to color
appropriately $V$ with $r(s-1)+1$ colors.
First, for each $i \in  \{0, \cdots, r-1\}$, we color the vertices
in independent sets $V_i(j)$ with color $j$,
for $j \in [\chi(G_i)]$.
This step, by the inductive assumption, can use exactly
$(r-1)(s-1)+1$ colors. Next, also by
the inductive assumption, we color properly
the graphs $H(j_0, \cdots, j_{r-1})$ just themselves
using $r(s-2)+1$ colors.
Note that this is $r$ less than the total number of
allowed colors. Let $k_i$, for $0 \le i < r$,
denote these extra colors.
We use colors $k_i$ to recolor some vertices
of the $H$-graphs, and thus obtain a proper coloring of
the entire $G(r,s)$ with $r(s-1)+1$ colors, as follows:
For each $r$-tuple $(j_0, \cdots, j_{r-1})$,
for every vertex $v \in V(H(j_0, \cdots, j_{r-1}))$
which has color $j_i$, recolor $v$ with a new color $k_i$,
for each $i \in  \{0, \cdots, r-1\}$.
\end{proof}

\bigskip
It is interesting to see that if the vertex of
$V_0$ is removed from $G(r,s)$, then
in the part (iii) of the proof of Theorem 3
the chromatic number
$\chi(G(r,s))$ drops to $m-1$, and thus
the arrowing of part (ii) would not hold.

We could improve a little the upper bound in Theorem 3
using an approach as in \cite{GS06},
but at the cost of significantly more complex 
construction, and hence we decided to not include it.

\bigskip
Recall that the chromatic Folkman graphs and numbers
we study here have just one additional constraint on
top of regular Folkman graphs and numbers, namely
the minimal possible chromatic number.
In general, in Theorem 3 and other places of this
paper, if $\chi$ is removed from an upper bound, then
it also holds for regular Folkman numbers. In particular, 
the special case of our Theorem 3 for two colors
but without considering chromatic numbers gives the bound
$F(2,s,s+1) \le 1 + s + sF(2,s-1,s)$.
The same follows from a construction by
Nenov \cite{Nenov85} using the corona product of graphs.
Bikov and Nenov \cite{BNc17} pointed out
that the constructions of \cite{Nenov85,KN05}
could be used as building blocks
in our proofs. Unfortunately, at the time of writing
the first version of this paper we were not aware of them.
Another upper bound construction of similar type
for vertex Folkman graphs was presented
by Xu et al. in \cite{GS06}.
Some special multicolor cases of Theorem~3
for regular Folkman numbers
(without considering chromatic numbers),
appear in Theorem~6 in \cite{LRU01} by
{\L}uczak, Ruci\'{n}ski and Urba\'{n}ski (2001).

\subsection{Bounds on $F_v^\chi$ for 2 colors}

A simple upper bound on $F^\chi(r,s,s+1)$ for $r=2$
based on Theorem 3 is $cs!$, for some positive constant $c$,
which is much larger than the upper bound on
$F(2,s,s+1)$ in Theorem~1. In this section we
look at some other 2-color cases: off-diagonal
$F_v^{\chi}(a, b; s+1)$ in Theorem 4 and
some special subcases of $F^\chi(s,s;s+1)$
in the sequel.

\medskip
\begin{theorem}
For any integers $a,b$ and $s$ such that $2 \le a,b \leq s$,
$F_v^{\chi}(a, b; s+1)$ exists and we have
$$F_v^{\chi}(a, b; s+1) \leq {{a+b-1} \over {2s-1}}F_v^{\chi}(s,s; s+1).
\eqno{(2)}$$
\end{theorem}

\begin{proof}
Suppose that $G \in \mathcal{F}_v^{\chi}(s,s; s+1)$,
and the order of $G$ is $F_v^{\chi}(s,s; s+1)$.
We clearly have $\chi(G)=2s-1$.
Write the set of vertices of $G$ as a partition
$$V(G) = \bigcup_{j=1} ^{2s-1} I_j,$$
where $I_j$'s are independent sets for
$j \in [2s-1]$, and $|I_{j_1}| \leq |I_{j_2}|$
for $j_1 < j_2$ and $j_1, j_2 \in [2s-1]$.
Let $G_i$ be the subgraph of $G$ induced by
$\bigcup_{j=1} ^i I_j$.
Note that this implies $\chi(G_i)=i$
for each $i$, since otherwise $\chi(G) < 2s-1$.

We claim that $G_{a+b-1} \in \mathcal{F}_v^{\chi}(a, b; s+1)$.
By the comments above we see that $\chi(G_{a+b-1})=a+b-1=m$
as required, and $cl(G_i)\le s$ holds by construction.
It remains to be shown that $G_{a+b-1} \rightarrow (a, b)^v$.
For a contradiction suppose that we have a red-blue
coloring of $V(G_{a+b-1})$ without any red $K_a$
and without any blue $K_b$. We can extend this coloring
to a full red-blue coloring of $V(G)$ by coloring
red all the vertices in $I_j$'s for
$a+b \le j \le s+b-1$, and coloring blue
all the vertices in $I_j$'s for
$s+b \le j \le 2s-1$. This coloring does not contain
any monochromatic $K_s$, which contradicts
the assumption that
$G \in \mathcal{F}_v^{\chi}(s,s; s+1)$.
Considering the non-decreasing orders of the sets $I_j$,
we have that
$$|V(G_i)|\le {{a+b-1} \over {2s-1}}F_v^{\chi}(s,s; s+1),$$
and thus the bound (2) and the theorem follow.
\end{proof}

\bigskip
Bollob\'{a}s and Thomason \cite{BB2006a} studied
the set-coloring of graphs, where an $r$-set-coloring
of a graph $G$ is defined as an assignment of
$r$ distinct colors to each vertex of $G$ so that the sets
of colors assigned to adjacent vertices are disjoint.
The set-coloring variants of the vertex- and edge
Folkman numbers were previously introduced and studied
by the first two authors of this paper jointly with
Wenfei Zhao and Zehui Shao \cite{setFolkman}.

\medskip
Let us denote the minimum number
of colors required to $r$-set-color any given graph $G$
by $\chi^{(r)}(G)$.
In 1979, Bollob\'{a}s and Thomason proved that
$\min \{ \chi ^{(r)}(G) \;|\; \chi(G) = t \}= t + 2r - 2$
\cite{BB2006a}.
We need a simple lemma using this result as follows.

\begin{lemma} \label{lemma5}
$\chi ^{(2s-1)}(C_{4s-1})= 4s-1$.
\end{lemma}

\begin{proof}
Using the above result
by Bollob\'{a}s and Thomason
on $\chi^{(r)}(G)$,
since $\chi(C_{4s-1}) = 3$, we clearly have
$\chi ^{(2s-1)}(C_{4s-1}) \geq 4s-1$.
On the other hand, it is easy to give
a proper $(2s-1)$-set-coloring witnessing
$\chi ^{(2s-1)}(C_{4s-1}) \leq 4s-1$.
We take both the vertices of the cycle and colors
to be in the set $\mathbb{Z}_{4s-1}$, and
assume that the edges of the cycle are $\{i,i+1\}$
for $i \in \mathbb{Z}_{4s-1}$,
all modulo $4s-1$. We assign the colors
$\{i(2s-1)+j\;|\;0 \le j \le 2s-2\}$ to
the vertex $i$ of $C_{4s-1}$, for each
$i \in \mathbb{Z}_{4s-1}$. One can easily see
that the sets of colors assigned to adjacent vertices
are disjoint.
\end{proof}

\medskip
The composition of simple graphs $G$ and $H$ is
denoted by $G[H]$, and it is defined as the graph
with vertex set $V(G) \times V(H)$, in which
vertex $(u,v)$ is adjacent to
$(u', v')$ if and only if either $u u' \in E(G)$
or $u = u'$ and $v v' \in E(H)$.
We will need another simple lemma involving
$G[H]$ and $\chi^{(r)}(G)$, namely:

\begin{lemma}
If $G$ and $H$ are graphs and $\chi(H) = r$, then
$\chi^{(r)}(G) = \chi(G[H])$.
\end{lemma}

\begin{proof}
Klav\v{z}ar proved that if $\chi(H) = r$, then
$\chi(G[H]) = \chi(G[K_r])$ \cite{K1996a}.
We can also easily see that
$\chi^{(r)}(G) =  \chi(G[K_r])$.
The lemma follows.
\end{proof}

\medskip
The bound $F_v(2s,2s; 2s+1) \leq 5F_v(s,s; s+1)$
was obtained constructively by Kolev
who used the composition of graphs $G[H]$
\cite{Nkolev2008}.
We will use a similar approach to obtain an upper
bound on the chromatic Folkman numbers of the form
$F_v^{\chi}(2s,2s; 2s+1)$.

\begin{theorem}
\label{th10} For any integer $s \geq 2$, we have
$$F_v^{\chi}(2s,2s; 2s+1) \leq (4s-1)F_v^{\chi}(s,s; s+1).$$
\end{theorem}

\begin{proof}
Let $H$ be any graph in $\mathcal{F}_v^{\chi}(s,s; s+1)$
of order $F_v^{\chi}(s,s; s+1)$, and thus $\chi(H) =2s-1$.
Observe that $C_{4s-1}[H] \rightarrow (2s,2s)^{v}$.
By Lemmas 5 and 6, we have
$\chi(C_{4s-1}[H]) = \chi ^{(2s-1)}(C_{4s-1})=4s-1$.
This in turn implies that
$C_{4s-1}[H] \in \mathcal{F}_v^\chi(2s,2s; 2s+1)$.
Finally, since the order of $C_{4s-1}[H]$ is equal to
$(4s-1)F_v^{\chi}(s,s; s+1)$ and clearly
$F_v^{\chi}(2s,2s; 2s+1) \leq |V(C_{4s-1}[H])|$,
this completes the proof.
\end{proof}

\medskip
\subsection{A conjecture}

\medskip
The classical Tur\'{a}n graph $T_ {n,r}$ is a complete multipartite graph
on $n$ vertices whose $r$ partite sets have sizes as equal as possible.
We say that a graph $G$ is {\em equitably $r$-colorable} if its vertices
can be properly colored with $r$-colors so that the sizes of any two
color classes differ by at most 1. Clearly, any $n$-vertex graph $G$
is equitably $r$-colorable if and only if $G$ is a subgraph of
$T_ {n,r}$. This motivates the formulation of the following
conjecture.

\bigskip
\noindent
{\bf Conjecture.}
{\it
For any integer $s \geq 2$, let $n = F_v(s,s; s+1)$. Then
there exists an $n$-vertex $K_{s+1}$-free subgraph $G$ of the Tur\'{a}n
graph $T_ {n, 2s-1}$, such that $G \rightarrow (s,s)^v$.}

\bigskip
Any subgraph of $T_ {n, 2s-1}$ has the
chromatic number upper bounded by $2s-1$. Therefore,
we can easily see that the Conjecture implies
the equality $F_v^{\chi}(s,s; s+1) = F_v(s,s; s+1)$,
in particular it would give a positive answer
to Problem 1 for two colors ($r=2$).
The Conjecture holds easily for $s=2$ and $n=5$ by
considering $C_5$. For $s=3$, we note that the graph
in Figure 1 of \cite{PRU99} witnessing
$F^{\chi}_v(3,3;4) = F_v(3,3;4) = 14$ is equitably 5-colorable,
and so it gives the second positive case. The Conjecture,
if proven true, would make the search for the upper
bound witnesses for $F_v(s,s; s+1)$ much easier,
including an approach using computer constructions.

Nenov studied several problems related to $F_v(r, 2, s+1)$,
for instance in \cite{Nenov2009a}.
Observe that using essentially only the definitions,
we can easily see that
$F_v(2s-2, 2, s+1) \leq F_v^{\chi}(s,s; s+1)$,
though we suspect that much better upper bound on
$F_v(2s-2, 2, s+1)$ is true.
Independently, $F_v(2s-2, 2, s+1)$ may be much
smaller than $F_v(s,s; s+1)$.
In another direction, one could use $C_5$, similarly
as we used $C_{4k-1}$ in Theorem 7, to study
the cases of $F_v^{\chi}(2s,2s,2s; 2s+1)$ and
more general diagonal and non-diagonal cases.
All of these problems seem interesting but difficult.

\medskip
\section{Minimum degree of graphs in $\mathcal{F}_v(s,s;s+1)$}
\label{Minimal}

\medskip
In this section we prove a theorem and then pose a problem
concerning lower bound on the minimum degree in some
minimal Folkman graphs. We consider only the case of
$F_v(s,s; s+1)$.

\begin{theorem}
For all integers $s \ge 3$, we have:

\noindent
{\rm (a)}
For every graph $G$, if $G \rightarrow (s,s)^v$ and
$G-u \not\rightarrow (s,s)^v$
for every vertex $u \in V(G)$,
then the minimum degree $\delta(G)$ satisfies
$\delta(G) \geq 2s-2$, and

\noindent
{\rm (b)}
There exists a $K_{s+1}$-free graph $G$ with minimum
degree $\delta(G) = 2s-2$, such that
$G \rightarrow (s,s)^v$ and
$G - u \not\rightarrow (s,s)^v$
for every vertex $u \in V(G)$.
\end{theorem}

\begin{proof}
(a)
For contradiction, let $u \in V(G)$ be any vertex
of degree at most $2s-3$. Assuming that
$G - u \not\rightarrow (s,s)^v$,
consider any $K_s$-free bipartition
$V_1 \cup V_2$ of the remaining vertices,
so that $V(G)=V_1 \cup V_2 \cup \{u\}$.
Without loss of generality we can also assume that
$|V_1 \cap N_G(u)| \le s-2$.
Color the vertices in $V_1 \cup \{u\}$ red and those in
$V_2$ blue, and note that this coloring has no
monochromatic $K_s$.
This contradicts $G \rightarrow (s,s)^v$,
and thus the part (a) follows.

\medskip
\noindent
(b)
Let $H$ be any $K_{s+1}$-free graph such that
$H \rightarrow (s,s)^v$. Assume further that $H$
is both vertex- and edge-minimal with respect to
arrowing, that is to say,
$H - u \not\rightarrow (s,s)^v$ for all vertices
$u \in V(H)$ and $H-e \not\rightarrow (s,s)^v$
for all edges $e \in E(H)$.
Fix some vertex $u_0 \in V(H)$, and
let $A$ be the set of all $(2s-2)$-element
subsets of $N_H(u_0)$ that contain two vertex-disjoint
$K_{s-1}$'s. Suppose that $A$ consists of $m$ sets,
i.e. $A = \{V_i \;|\; 1 \leq i \leq m \}$.
Note that vertex-minimality of $H$ implies
that $A$ is nonempty, so $m \ge 1$.
First, we extend graph $H$ to $H'$ by
adding new vertices $\{u_i \;|\; 1 \leq i \leq m \}$
and edges $\{u_i v \;| \; v \in V_i \}$ for
all $i \in [m]$.
Next, we delete vertex $u_0$ from $H'$.
Observe that for all $i \in [m]$ the degree
of vertex $u_i$ in $H'$ is equal to $2s-2$.

We claim that $H' \rightarrow (s,s)^v$. Suppose
to the contrary, namely that there exists a red-blue
coloring $C'$ of $V(H')$ without any monochromatic $K_s$.
Thus, the restriction of $C'$ to a coloring
of the vertices of $H-u_0$, say $C$, is a witness of
$H - u_0 \not\rightarrow (s,s)^v$. Since
$H \rightarrow (s,s)^v$, $C$ must contain
two vertex-disjoint monochromatic $K_{s-1}$'s,
furthermore they must be in different colors,
and both are contained in $N_H(u_0)$. This however
contradicts the properties of $C'$ following
from the construction of $H'$.
Hence $H' \rightarrow (s,s)^v$.

We will define the final graph $G$ satisfying (b) to
be an induced subgraph of $H'$  on the vertex set
of the form $V(H') \setminus B$, where
$B \subset \{ u_i \;|\; 1 \le i \le m \}$.
We choose $B$ so that its vertex indices
form a maximal subset of $[m]$
still giving $G \rightarrow (s,s)^v$.
The properties of $H'$ stated above guarantee
that such $B$ must be a proper subset of $[m]$.
This can be seen, since if we delete
$\{u_i \:|\: 1 \leq i \leq m \}$,
we obtain a graph isomorphic to $H - u_0$,
which does not arrow $(s,s)^v$.
This way, we can obtain a minimal graph
$G$ which arrows $(s,s)^v$ with at least
one vertex $u_i$ of degree $2s-2$.
Together with part (a), this completes
the proof of (b).
\end{proof}

\medskip
Bikov and Nenov suggested \cite{BNc17} that in our proof
of Theorem 8(b) we could use graphs $M_k$, for odd $k$,
defined and studied by Nenov (\cite{Nenov85}, page 351).

\medskip
Finally, we pose the following question.

\medskip
\noindent
{\bf Problem 2.} {\em
For which integers $n, s \ge 3$
does there exist a $(2s-2)$-regular $K_{s+1}$-free graph $G$
on $n$ vertices such that $G \rightarrow (s,s)^v?$
}

\medskip
It seems that even the case of $s=3$ is not obvious.
If the answer to this problem for each $s$ is YES for
at least some $n$, then the chromatic number of such
a graph must be equal to $2s-1$. Thus, it could give
another proof of the existence of $F_v^{\chi}(s,s; s+1)$.
Of course, the order $n$ of such a graph $G$ may be much
larger than $F_v^{\chi}(s,s; s+1)$.

\smallskip
\section*{Acknowledgements}

We are grateful to Nedyalko Nenov and Aleksandar Bikov
for their discussions and comments on an earlier version
of this paper.
We also would like to thank the anonymous reviewer for
detailed suggestions of how to revise the paper, which
led to a very significant improvement of the presentation.

\smallskip

\end{document}